\documentclass[11pt]{article}
\ifdefined\pdfminorversion\pdfminorversion=7\fi

\usepackage[a4paper,margin=24mm]{geometry}
\usepackage[T1]{fontenc}
\usepackage[utf8]{inputenc}
\usepackage{lmodern}
\usepackage{microtype}
\usepackage{amsmath,amssymb,amsthm,mathtools}
\usepackage{enumitem}
\usepackage{graphicx}
\usepackage{tikz}
\usepackage{tikz-cd}
\usetikzlibrary{arrows.meta,positioning,calc}
\usepackage[hidelinks]{hyperref}
\usepackage[nameinlink,noabbrev]{cleveref}
\hypersetup{
 pdftitle={Convergence of Nonmonotone Schemes to Viscosity Solutions via Local Comparison},
 pdfauthor={Alessio Basti and Fabio Camilli}
}

\allowdisplaybreaks
\setlist[enumerate]{leftmargin=*,itemsep=0.25em,topsep=0.4em}

\newtheorem{theorem}{Theorem}[section]
\newtheorem{proposition}[theorem]{Proposition}
\newtheorem{lemma}[theorem]{Lemma}
\newtheorem{corollary}[theorem]{Corollary}
\newtheorem{definition}[theorem]{Definition}
\theoremstyle{remark}
\newtheorem{remark}[theorem]{Remark}

\newcommand{\halfup}{\overline u}
\newcommand{\halfdown}{\underline u}

\newcommand{\Qh}{\mathcal Q_h}

\newcommand{\QTwoFinestN}{512}
\newcommand{\QTwoLastOrder}{1.98}
\newcommand{\QTwoFirstError}{$4.73\times10^{-3}$}
\newcommand{\QTwoLastError}{$4.76\times10^{-6}$}
\newcommand{\QTwoMaximumResidual}{$6.02\times10^{-10}$}

\newcommand{\QTwoStrictFractions}{17.22\%, 23.87\%, 48.73\%, 10.18\%}
\newcommand{\QTwoOrderGap}{3.58}
\newcommand{\QTwoOrderGapFraction}{18.98\%}

\numberwithin{equation}{section}

\title{Convergence of Nonmonotone Schemes to Viscosity Solutions via Local Comparison}
\author{Alessio Basti\footnote{Engineering and Geology Department, ``G. d'Annunzio'' University Chieti-Pescara,
		viale Pindaro 42, 65127 Pescara (Italy), { \textsc alessio.basti@unich.it, fabio.camilli@unich.it.}} \and
 Fabio Camilli$^*$}

\begin{document}
\maketitle
\begin{abstract}

We develop a convergence criterion for nonmonotone approximation schemes
for Hamilton--Jacobi and fully nonlinear second-order equations.
Under stability, consistency and suitable boundary assumptions, a local
comparison hypothesis replaces monotonicity and yields convergence to
viscosity solutions. For a semi-Lagrangian CWENO scheme for one-dimensional
Hamilton--Jacobi equations, we establish convergence for bounded uniformly
continuous initial data under $h=o(\tau)$, relaxing the previously required
mesh condition $h=O(\tau^2)$. We also prove convergence of a nonmonotone
one-dimensional $Q^2$ scheme for Isaacs equations whose operator cannot
be written as a monotone operator plus a perturbation vanishing uniformly
on bounded grid functions.

\end{abstract}



\section{Introduction}

The convergence theorem of Barles and Souganidis
\cite{BarlesSouganidis1991} is a basic tool in the numerical analysis of
Hamilton--Jacobi equations and fully nonlinear second-order equations.
Stability, consistency, and monotonicity, together with comparison for the
limiting equation, yield convergence to the viscosity solution.
The order-preserving property is difficult to retain in higher-order
polynomial reconstructions and finite element stencils.
Positive-type discretizations of second-order HJB and fully nonlinear
equations were studied in \cite{CamilliFalcone1995,KuoTrudinger1992}.
Oberman \cite{Oberman2006} discussed monotonicity criteria and showed
that stability alone need not ensure convergence; see also
\cite{BarlesJakobsen2002} for convergence rates of monotone HJB schemes.

We write the discrete equation as $S(\rho,x,r,V)=0$, where $\rho$
denotes the discretization parameters, $x$ is a grid node, $r$ is the
value at that node, and $V$ contains the grid values entering the scheme.
With our sign convention, Barles--Souganidis monotonicity reads
\begin{equation}
	V\ge W\quad\Longrightarrow\quad
	S(\rho,x,r,V)\le S(\rho,x,r,W).
	\label{eq:BS}
\end{equation}
Several approaches have been developed to deal with schemes that do not
satisfy this condition, including filtered and nearly monotone schemes
\cite{FroeseOberman2013,BokanowskiFalconeSahu2016},
$\varepsilon$-monotonicity and related constructions
\cite{BokanowskiEtAl2015,Warin2019}, filtered approximations for
second-order HJB equations \cite{BokanowskiPicarelliReisinger2018},
and methods based on numerical moments and
g-monotonicity \cite{FengLewis2021,FengLewisWard2022}.
Another relaxation of monotonicity, based on comparison
with smooth test functions and allowing for a different
consistent comparison operator, was used in
\cite{CarliniFalconeFerretti2010,AchdouFalcone2012}.
In a different direction, Jensen and Smears \cite{JensenSmears2013}
used elliptic projections to establish convergence of monotone
finite element schemes for HJB equations that do not satisfy the usual consistency condition;
see also \cite{SalgadoZhang2019} for a two-scale approximation
of Isaacs equations.

Our starting point is that classical monotonicity is used when a
numerical solution is compared with a smooth test function near a
half-relaxed-limit contact point. We therefore require only the
comparison needed in this situation, on finite grid patches and their
stencil boundaries. If the numerical family crosses a shifted test
function inside a patch while remaining ordered on its boundary, a
suitable discrete inequality must hold at some interior node, not
necessarily at a maximal crossing. Intermediate shifts and vanishing
comparison errors are allowed. For earlier local comparison arguments,
see, e.g., \cite{NochettoNtogkasZhang2019}.

We incorporate this local
comparison into a convergence criterion based on half-relaxed limits.
The criterion retains the usual stability and consistency requirements,
but allows comparison to be verified for an auxiliary numerical family
and a local auxiliary operator. The auxiliary family must remain close
to the original solution and satisfy the auxiliary equation up to a
vanishing error; the auxiliary operator must be consistent with the
same limiting equation. This makes it possible to recover the viscosity
inequalities without assuming monotonicity of the original scheme.
The usefulness of the criterion depends, of course, on whether its
local comparison hypothesis can be checked in concrete examples.

The applications illustrate two sufficient conditions. For the semi-Lagrangian CWENO scheme of
\cite{CarliniFerrettiPredaSemplice2025}, we relax the mesh condition
from $h=O(\tau^2)$ to $h=o(\tau)$, allowing bounded uniformly continuous
initial data. The proof uses the \emph{monotone transfer mechanism}, which gives the local comparison by
\cref{prop:local-transfer}. For the nonmonotone one-dimensional $Q^2$ Isaacs scheme, we prove
convergence outside the class of monotone operators with perturbations
vanishing uniformly on bounded grid functions. Here, local comparison
follows from a \emph{strict patchwise comparison} on
aligned grid intervals, by \cref{prop:all-routes}.

The paper is organized as follows. Section~\ref{sec:criterion}
introduces the local comparison hypothesis and proves the convergence
criterion. Section~\ref{sec:verification} gives sufficient conditions
based on monotone transfer and patchwise comparison.
Section~\ref{sec:cweno} establishes the CWENO convergence
result, while Section~\ref{sec:q2} develops the
one-dimensional $Q^2$ Isaacs application. Numerical tests are reported
in Section~\ref{sec:numerics}, and Appendix~\ref{app:technical-proofs}
contains the technical auxiliary proofs.

\section{A convergence criterion beyond monotonicity}
\label{sec:criterion}

Let $\Omega\subset\mathbb R^d$ be a bounded open set. We consider the fully nonlinear second order equation
\begin{equation}
	\begin{cases}
		F(D^2u(x),Du(x),u(x),x)=0, & x\in\Omega,\\
		u=g, & x\in\partial\Omega,
	\end{cases}
	\label{eq:pde}
\end{equation}
where $F:\mathbb S^d\times\mathbb R^d\times\mathbb R\times\overline\Omega
\longrightarrow\mathbb R$ is a possibly discontinuous degenerate elliptic operator. Here
$\mathbb S^d$ denotes the space of real symmetric $d\times d$ matrices.
When $F$ is discontinuous, we denote by $F_*$ and $F^*$ its lower and
upper semicontinuous envelopes, respectively.
\begin{definition}
 An upper
semicontinuous function $u$ is a viscosity subsolution in $\Omega$ of \eqref{eq:pde} if,
whenever $\phi\in C^2(\Omega)$ and $u-\phi$ has a local maximum at
$x\in\Omega$,
\[
F_*(D^2\phi(x),D\phi(x),u(x),x)\le 0.
\]
Similarly, a lower semicontinuous function $u$ is a viscosity
supersolution of  \eqref{eq:pde} if, whenever $u-\phi$ has a local minimum at $x$,
\[
F^*(D^2\phi(x),D\phi(x),u(x),x)\ge 0.
\]
A continuous function which is both a subsolution and a supersolution
is a viscosity solution.
\end{definition}
See \cite{CrandallIshiiLions1992} for viscosity solutions.
The theorem below concerns the interior; boundary conditions and
convergence are treated in \cref{rem:boundary,cor:convergence}.

We now describe the discrete approximation. Let
$\mathcal R\subset(0,\infty)$ be a set of discretization parameters
having $0$ as an accumulation point. For each $\rho\in\mathcal R$, let
$\mathcal G_\rho\subset\overline\Omega$ be a finite grid satisfying
\begin{equation}
	\sup_{x\in\overline\Omega}
	\operatorname{dist}_\star(x,\mathcal G_\rho)\longrightarrow0
	\qquad\text{as }\rho\to0,
	\label{eq:density}
\end{equation}
where $|\cdot|_\star$ is a fixed norm on $\mathbb R^d$. All limits below
are understood along $\rho\in\mathcal R$.
A numerical scheme is written in the form
\begin{equation}
	S(\rho,x,U^\rho(x),U^\rho)=0,
	\qquad x\in\mathcal G_\rho\cap\Omega,
\end{equation}
where $U^\rho:\mathcal G_\rho\to\mathbb R$ is a grid function.
Here $\rho$ denotes the discretization, $x$ the grid node,
$r=U^\rho(x)$ the central value, and $V=U^\rho$ the grid values.
Write $\mathcal S_\rho[V](x):=S(\rho,x,V(x),V)$.
We assume locality: for every
$x\in\mathcal G_\rho\cap\Omega$ there is a finite set
\[
\operatorname{St}_\rho(x)\subset\mathcal G_\rho,
\qquad x\in\operatorname{St}_\rho(x),
\]
called the stencil of the scheme at $x$, such that
$\mathcal S_\rho[V](x)$ depends on $V$ only through its values on
$\operatorname{St}_\rho(x)$. 
We also require the stencil diameter to vanish locally: for every
compact set $K\Subset\Omega$,
\begin{equation}
	\sup_{x\in\mathcal G_\rho\cap K}
	\sup_{y\in\operatorname{St}_\rho(x)}
	|y-x|_\star
	\longrightarrow0.
	\label{eq:locality}
\end{equation}
We retain stability and consistency and replace monotonicity by local comparison.

\medskip
\noindent\textbf{(H1) Stability.}
For all sufficiently small $\rho$, the discrete problem admits a grid function
$u^\rho:\mathcal G_\rho\to\mathbb R$ satisfying the prescribed discrete
boundary conditions and
\[
\mathcal S_\rho[u^\rho](x)=0,
\qquad x\in\mathcal G_\rho\cap\Omega.
\]
Moreover, the family of discrete solutions is uniformly bounded:
\[
\sup_{\rho}
\|u^\rho\|_{\ell^\infty(\mathcal G_\rho)}<\infty.
\]

\medskip
\noindent\textbf{(H2) Consistency.}
For every $x\in\Omega$ and every $\phi\in C_b^2(\mathbb R^d)$,
\begin{align}
	\limsup_{\substack{\rho\downarrow0,\ y\to x,\ \xi\to0\\
			y\in\mathcal G_\rho\cap\Omega}}
	\mathcal S_\rho[\phi+\xi](y)
	&\le
	F^*(D^2\phi(x),D\phi(x),\phi(x),x),
	\notag\\
	\liminf_{\substack{\rho\downarrow0,\ y\to x,\ \xi\to0\\
			y\in\mathcal G_\rho\cap\Omega}}
	\mathcal S_\rho[\phi+\xi](y)
	&\ge
	F_*(D^2\phi(x),D\phi(x),\phi(x),x).
	\label{eq:consistency}
\end{align}
Here, as throughout the paper, a smooth function is identified with its
restriction to the grid.

The operator is normalized as in \eqref{eq:consistency}: it corresponds
to $S_{\mathrm{BS}}/\rho$ in \cite{BarlesSouganidis1991}.
This normalization matters when discrete inequalities include errors
that must vanish with the mesh.
\subsection{Local comparison}

A \emph{grid patch} is a nonempty finite set
$D\subseteq\mathcal G_\rho\cap\Omega$.
Because the stencils of nodes in $D$ may extend outside it,
we define its stencil boundary by
\[
\partial_\rho D
:=
\left(
\bigcup_{x\in D}\operatorname{St}_\rho(x)
\right)\setminus D,
\].
Thus $\partial_\rho D$ contains precisely the nodes outside $D$
needed to evaluate the scheme inside it. For each $\rho$ we select
a family $\mathfrak D_\rho$ of such patches.

\begin{definition}[Localizing family of patches]
	We say that $\mathfrak D=(\mathfrak D_\rho)_\rho$
	is a \emph{localizing family of patches} if, for every $x_0\in\Omega$ and every pair of radii
	$0<r<R$ such that $\overline B_R(x_0)\Subset\Omega$,
	there exists, for all sufficiently small $\rho$, a grid patch
	$D\in\mathfrak D_\rho$ satisfying
	\begin{equation}
		\mathcal G_\rho\cap B_r(x_0)
		\subseteq D
		\subseteq
		\mathcal G_\rho\cap B_R(x_0),
		\qquad
		\partial_\rho D
		\subseteq
		\mathcal G_\rho\cap
		\bigl(B_R(x_0)\setminus B_r(x_0)\bigr).
		\label{eq:basis}
	\end{equation}
\end{definition}

\begin{figure}[htbp]
	\centering
	\includegraphics[width=0.90\linewidth]{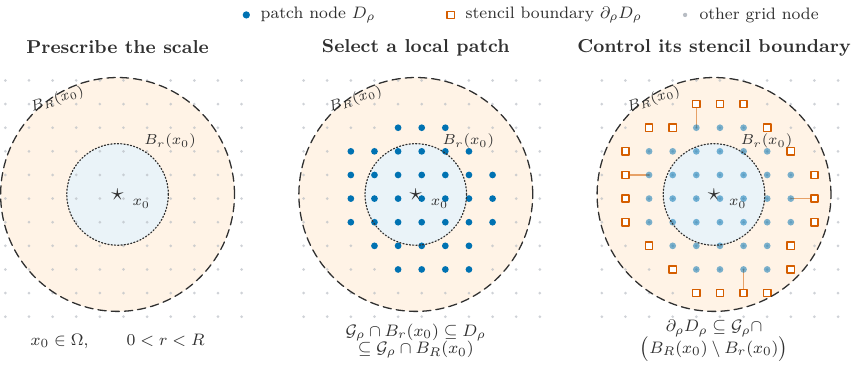}
	\caption{A localizing patch $D$ contains the inner grid ball, while its stencil boundary $\partial_\rho D$ lies in $B_R(x_0)\setminus B_r(x_0)$.}
	\label{fig:patch-basis}
\end{figure}

Condition \eqref{eq:basis} keeps the patch within $B_R(x_0)$,
contains $B_r(x_0)\cap\mathcal G_\rho$, and separates its stencil
boundary from $x_0$; see \cref{fig:patch-basis}.

The balls in \eqref{eq:basis} are taken in the fixed norm
$|\cdot|_\star$. In applications, the family must also support a required comparison
argument that will be subsequently described; the admissible patch geometry is therefore part of each application, as illustrated by the aligned intervals in \cref{sec:q2}.
A general geometric construction useful for the applications may be to choose an
intermediate radius $r<s<R$ and set
\(
 D=\mathcal G_\rho\cap B_s(x_0).
\)
By \eqref{eq:locality}, for sufficiently small $\rho$ the stencil radius
on $\overline B_R(x_0)$ is smaller than $R-s$. Every stencil node used
in $D$ then lies in $B_R(x_0)$. Since $D$ contains every grid node of
$B_s(x_0)$, its stencil boundary lies in $B_R(x_0)\setminus B_s(x_0)$,
which gives \eqref{eq:basis}. Grid density ensures that $D$ is nonempty
for small $\rho$, and the finiteness of the grid makes it a finite patch.
Taking all such grid balls therefore provides a localizing family.

Let $\mathcal A_\rho$ be a local discrete operator and $Z^\rho$
a numerical family. Denote their stencils, stencil boundaries, and
localizing patches by $\operatorname{St}^{\mathcal A}_\rho(x)$,
$\partial^{\mathcal A}_\rho D$, and $\mathfrak D^{\mathcal A}$.
The auxiliary stencils may differ from those of the original scheme.

Discrete comparison (DC) requires an operator inequality when a
shifted smooth test function and $Z^\rho$ are ordered on the stencil
boundary but cross inside the patch.

\begin{definition}[DC]
\label[definition]{def:dc}
We say that $(\mathcal A_\rho,Z^\rho)$ satisfies
\textnormal{(DC)} if there exists a localizing family of patches
$\mathfrak D^{\mathcal A}$ with the following property. For every compact $K\Subset\Omega$ and every test function
$\phi\in C_b^2(\mathbb R^d)$, there exist $\rho_0,c_0>0$ such that,
for every $0<\rho<\rho_0$, every $|c|<c_0$, and every patch
$D\in\mathfrak D^{\mathcal A}_\rho$ satisfying $D\cup\partial^{\mathcal A}_\rho D\subset K$,
we require:
\begin{enumerate}[
    label=\textnormal{(\roman*)},
    leftmargin=*,
    itemsep=0.4em
]
\item \emph{Upper crossing.}
If
\[
    Z^\rho\le\phi+c
    \quad\text{on }\partial^{\mathcal A}_\rho D,
    \qquad
    m:=\max_D(Z^\rho-\phi-c)>0,
\]
then there exist $z\in D$ and $\theta\in[0,m]$ such that
\begin{equation}
    \mathcal A_\rho[Z^\rho](z)
    \ge
    \mathcal A_\rho[\phi+c+\theta](z).
    \label{eq:dc-upper}
\end{equation}

\item \emph{Lower crossing.}
If
\[
    \phi+c\le Z^\rho
    \quad\text{on }\partial^{\mathcal A}_\rho D,
    \qquad
    m:=\max_D(\phi+c-Z^\rho)>0,
\]
then there exist $z\in D$ and $\theta\in[0,m]$ such that
\begin{equation}
    \mathcal A_\rho[Z^\rho](z)
    \le
    \mathcal A_\rho[\phi+c-\theta](z).
    \label{eq:dc-lower}
\end{equation}
\end{enumerate}
\end{definition}

The two conditions correspond to crossings from below and from above;
$\theta$ allows an intermediate vertical shift. 

For a Barles--Souganidis monotone operator, \textnormal{(DC)} holds.
In the upper-crossing case, let $z$ maximize
$Z^\rho-\phi-c$ in $D$ and take $\theta=m$. The two grid functions have
the same value at $z$, while
$Z^\rho\le\phi+c+m$ on the stencil of $z$; monotonicity gives
\eqref{eq:dc-upper}. The lower-crossing case is analogous. Thus
\textnormal{(DC)} retains precisely the local comparison step used in the
classical proof, but requires it only for the selected grid family and
smooth test functions, rather than for arbitrary ordered grid functions.

For a grid function $V$ we write \(\|V\|_K:=\|V\|_{\ell^\infty(K\cap\mathcal G_\rho)},\) with the convention that the norm is zero when
$K\cap\mathcal G_\rho$ is empty.

We replace the monotonicity assumption of the Barles--Souganidis
framework by the following hypothesis:

\medskip
\noindent\textbf{(H3) Local comparison.}
There exist a local discrete operator \(\mathcal M_\rho[V](x)=M(\rho,x,V(x),V)\) and grid functions $v^\rho:\mathcal G_\rho\to\mathbb R$ such that, for
every compact $K\Subset\Omega$, the following conditions hold.
\begin{enumerate}[label=\textnormal{(\alph*)}]
\item \emph{Consistency.} \textnormal{(H2)} applies to $\mathcal M_\rho$.

\item \leftprotrusion{\emph{Approximation and vanishing residual.}}
\begin{gather}
 \|v^\rho-u^\rho\|_K\longrightarrow0,
 \label{eq:auxiliary-close}\\
 \eta_\rho(K):=\|\mathcal M_\rho[v^\rho]\|_K
 \longrightarrow0.
 \label{eq:auxiliary-residual}
\end{gather}

\item \emph{Discrete comparison.}
The pair $(\mathcal M_\rho,v^\rho)$ satisfies \textnormal{(DC)} on a
localizing family of patches associated with the stencils of $\mathcal M_\rho$.
\end{enumerate}

Condition (a) identifies the limiting equation; (b) gives the same
interior half-relaxed limits and a vanishing auxiliary equation error;
(c) provides the order comparison needed for smooth tests.

For a monotone scheme one takes $\mathcal M_\rho=\mathcal S_\rho$
and $v^\rho=u^\rho$; (DC) also allows a different comparison node
and an intermediate shift $\theta$.

\begin{remark}[The auxiliary residual and the shift $\theta$]

The defect $\eta_\rho(K)$ of \eqref{eq:auxiliary-residual} is essential:
in \cref{sec:cweno}, the numerical solution of the original scheme
satisfies the monotone auxiliary equation only up to a locally uniform
$o(1)$ residual.
The interval $[0,m]$ for the shift $\theta$ allows a single formulation to
cover both verification mechanisms: monotonicity produces
\eqref{eq:dc-upper}--\eqref{eq:dc-lower} with $\theta=m$, whereas patchwise
comparison produces them with $\theta=0$.
\end{remark}

\medskip
Define the half-relaxed limits of the original numerical family by
\begin{equation}
 \halfup(x)
 :=\limsup_{\substack{\rho\downarrow0,\ y\to x\\y\in\mathcal G_\rho}}
 u^\rho(y),
 \qquad
 \halfdown(x)
 :=\liminf_{\substack{\rho\downarrow0,\ y\to x\\y\in\mathcal G_\rho}}
 u^\rho(y).
\end{equation}

\begin{theorem}[Interior half-relaxed limits]
\label[theorem]{thm:interior}
Assume the grid-density and locality conditions
\eqref{eq:density}--\eqref{eq:locality}, together with
\textnormal{(H1)-(H3)}. Then $\halfup$ is a viscosity subsolution and
$\halfdown$ is a viscosity supersolution of
\[
 F(D^2u,Du,u,x)=0
 \qquad\text{in }\Omega.
\]
\end{theorem}

\begin{proof}
By \eqref{eq:auxiliary-close}, $u^\rho$ and $v^\rho$ have the same upper
and lower half-relaxed limits at every interior point. We may therefore
prove the viscosity inequalities using the auxiliary family. As usual, it
is enough to consider strict contacts. A local $C^2$ test function can be
modified away from the contact point so as to belong to
$C_b^2(\mathbb R^d)$.
Let $\phi\in C_b^2(\mathbb R^d)$ touch $\halfup$ strictly from above at
$x_0\in\Omega$, with $\phi(x_0)=\halfup(x_0)$. We prove
\[
 F_*(D^2\phi(x_0),D\phi(x_0),\phi(x_0),x_0)\le0.
\]
Suppose, to the contrary, that this quantity is positive. By lower
semicontinuity, and by shrinking the neighbourhood of the strict contact,
we may choose $0<r<R$, $\sigma>0$, and
$K=\overline B_R(x_0)\Subset\Omega$ such that
\begin{equation}
 F_*(D^2\phi(x),D\phi(x),\phi(x),x)\ge2\sigma
 \qquad(x\in K),
 \label{eq:positive-jet}
\end{equation}
and
\[
 \halfup-\phi<0
 \qquad\text{on }K\setminus B_r(x_0).
\]
Compactness of the annulus and the definition of the upper half-relaxed
limit yield $\gamma>0$ such that, for all sufficiently small $\rho$,
\begin{equation}
 v^\rho-\phi\le-\gamma
 \qquad\text{on }\mathcal G_\rho\cap
       (K\setminus B_r(x_0)).
 \label{eq:annulus-upper}
\end{equation}

Choose $\rho_n\downarrow0$ and $x_n\to x_0$ such that
\[
 d_n:=v^{\rho_n}(x_n)-\phi(x_n)\longrightarrow0,
\]
and set \(c_n:=d_n-\frac1n.\) Since $c_n\to0$, we may restrict to indices $n$ so large that
$|c_n|<\min\{\gamma,c_0\}$, with $c_0$ as in \cref{def:dc}; the shift $c_n$
is then admissible in \textnormal{(DC)} and does not destroy the ordering on
the annulus. For each such $n$, choose from the localizing family a patch
$D_n$ associated with the inner ball $B_r(x_0)$ and the outer ball
$B_R(x_0)$. For large $n$, $x_n\in D_n$, and \eqref{eq:annulus-upper} gives
\[
 v^{\rho_n}\le\phi+c_n
 \qquad\text{on }\partial^{\mathcal M}_{\rho_n}D_n.
\]
At $x_n$ the reverse inequality is strict, since
$v^{\rho_n}(x_n)-\phi(x_n)-c_n=1/n$. Hence
\[
 m_n:=\max_{D_n}(v^{\rho_n}-\phi-c_n)>0.
\]
Moreover, $m_n\to0$. Otherwise, a subsequence of maximizing nodes would
converge to a point of $K$ at which $\halfup-\phi>0$, contradicting the
choice of the test function.

Apply \eqref{eq:dc-upper} with the shift $c_n$. There exist $z_n\in D_n$
and $\theta_n\in[0,m_n]$ such that
\begin{equation}
 \mathcal M_{\rho_n}[\phi+c_n+\theta_n](z_n)
 \le \mathcal M_{\rho_n}[v^{\rho_n}](z_n)
 \le \eta_{\rho_n}(K)
 \longrightarrow0.
 \label{eq:upper-total-defect}
\end{equation}
Here $c_n\to0$ and $0\le\theta_n\le m_n\to0$. After extracting a
subsequence, $z_n\to z\in K$. The lower relaxed consistency of
$\mathcal M_\rho$ gives
\[
 F_*(D^2\phi(z),D\phi(z),\phi(z),z)
 \le
 \liminf_{n\to\infty}
 \mathcal M_{\rho_n}[\phi+c_n+\theta_n](z_n),
\]
whereas \eqref{eq:upper-total-defect} makes the right-hand side
nonpositive. This contradicts \eqref{eq:positive-jet}. Thus $\halfup$ is
a viscosity subsolution.

The supersolution argument is symmetric. Let $\phi$ touch $\halfdown$
strictly from below at $x_0$ and assume, by contradiction,
\[
F^*(D^2\phi(x_0),D\phi(x_0),\phi(x_0),x_0)<0.
\]
Choose $0<r<R$, $\sigma>0$, and
$K=\overline B_R(x_0)\Subset\Omega$ so that
\[
F^*(D^2\phi,D\phi,\phi,x)\le-2\sigma
\qquad\text{on }K,
\]
and, for small $\rho$,
\[
v^\rho-\phi\ge\gamma
\qquad\text{on }\mathcal G_\rho\cap(K\setminus B_r(x_0)).
\]

Take $\rho_n\downarrow0$ and $x_n\to x_0$ such that
$d_n:=v^{\rho_n}(x_n)-\phi(x_n)\to0$, and set
\[
c_n=d_n+\frac1n.
\]
For the corresponding patches $D_n$, one has
\[
\phi+c_n\le v^{\rho_n}
\quad\text{on }\partial^{\mathcal M}_{\rho_n}D_n,
\qquad
m_n:=\max_{D_n}(\phi+c_n-v^{\rho_n})>0,
\]
with $m_n\to0$ by the same compactness argument as above. Hence
\eqref{eq:dc-lower} gives points $z_n\in D_n$ and
$0\le\theta_n\le m_n$ such that
\begin{equation}
	\mathcal M_{\rho_n}[\phi+c_n-\theta_n](z_n)
	\ge
	-\eta_{\rho_n}(K)
	\longrightarrow0.
\end{equation}
Passing to a subsequence and using the upper relaxed consistency of
$\mathcal M_\rho$ contradicts the strict negativity of $F^*$ on $K$.
Therefore $\halfdown$ is a viscosity supersolution.

This proves the viscosity subsolution and supersolution properties of
$\halfup$ and $\halfdown$ in $\Omega$.
\end{proof}
\begin{remark}[Boundary conditions]
	\label[remark]{rem:boundary}
	The preceding theorem is intentionally an interior result. Its proof uses
	only patches compactly contained in $\Omega$ and therefore does not address
	the behavior of the half-relaxed limits on $\partial\Omega$.
	
	To obtain convergence for a boundary-value problem associated with
	\eqref{eq:pde}, one must additionally verify that $\halfup$ and
	$\halfdown$ satisfy the prescribed boundary condition in a form compatible
	with the comparison principle. This boundary analysis depends on the
	equation, the geometry of the domain, and the discretization of the boundary
	data. For examples where the Dirichlet condition must be understood in the
generalized viscosity sense, and for the corresponding comparison
principles, see \cite{BarlesDaLio2004}. As in the Barles--Souganidis framework, it is logically separate from
	the local crossing mechanism used here to replace monotonicity in the
	interior.
\end{remark}

\begin{corollary}[Convergence]
	\label[corollary]{cor:convergence}
	\label[corollary]{thm:main}
	Assume the hypotheses of \cref{thm:interior}. Suppose in addition that
	$\halfup$ and $\halfdown$ are, respectively, a viscosity subsolution and a
	viscosity supersolution of the full boundary-value problem associated with
	\eqref{eq:pde}, in a class for which comparison holds. Then
	\[
	\halfup=\halfdown=:u
	\qquad\text{on }\overline\Omega,
	\]
	$u\in C(\overline\Omega)$, and
	\begin{equation}
		\max_{z\in\mathcal G_\rho}
		|u^\rho(z)-u(z)|
		\longrightarrow0.
		\label{eq:convergence}
	\end{equation}
\end{corollary}

\begin{proof}
	Comparison gives $\halfup\le\halfdown$ on $\overline\Omega$, while by
	definition $\halfdown\le\halfup$. Hence
	$\halfup=\halfdown=:u$, and $u$ is continuous since it is both upper and
	lower semicontinuous.
	
	To prove uniform convergence, suppose that \eqref{eq:convergence} fails.
	Then, for some $\varepsilon>0$, there are $\rho_n\downarrow0$ and
	$z_n\in\mathcal G_{\rho_n}$ such that
	\[
	|u^{\rho_n}(z_n)-u(z_n)|\ge\varepsilon.
	\]
	After extracting a subsequence, $z_n\to z\in\overline\Omega$. Stability
	and continuity of $u$, together with the definitions of the half-relaxed
	limits, imply
	\[
	\halfdown(z)
	\le
	\liminf_n u^{\rho_n}(z_n)
	\le
	\limsup_n u^{\rho_n}(z_n)
	\le
	\halfup(z)
	=
	u(z),
	\]
	and therefore $u^{\rho_n}(z_n)\to u(z)$, contradicting the preceding
	inequality.
\end{proof}

\begin{remark}[Evolutionary problems]
	
	The argument also applies to space--time grids: consistency then includes
the time derivative in $u_t+F(D_x^2u,D_xu,u,t,x)=0$, while the initial
condition is handled separately through comparison for the Cauchy
problem. Periodic grids need no lateral boundary analysis; see the
CWENO application in \cref{sec:cweno}.
\end{remark}
\begin{remark}[Generalized monotonicity]
	Condition (DC) can allow a different consistent operator on smooth
test functions and an additive error $\omega_\rho\ge0$: subtract
$\omega_\rho$ in the upper-crossing inequality and add it in the
lower-crossing inequality. The error must vanish jointly with the
mesh and constant shift, uniformly over admissible patches and
comparison points. The argument of \cref{thm:interior} then uses
consistency of the comparison operator. This is related to generalized
monotonicity in \cite{CarliniFalconeFerretti2010,AchdouFalcone2012}.
\end{remark}

\begin{remark}[Locally finite grids and unbounded domains]
	
	The interior proof also works on unbounded domains with locally
finite grids and finite patches, provided density, locality, and
patch conditions hold on compact sets and (H1) supplies uniform
stability. This includes the Cauchy problem in \cref{sec:cweno}.
\end{remark}

\begin{corollary}[Cauchy problems]
\label[corollary]{cor:cauchy}
Let $Q=(0,\infty)\times\mathbb R^d$ and let
$\mathcal G_\rho\subset\overline Q$ be locally finite grids,
asymptotically dense on compact subsets of $\overline Q$.
Assume that $u^\rho$ is uniformly bounded on $[0,T]\times\mathbb R^d$ for every $T>0$ and that \textnormal{(H3)} holds locally in space--time,
with consistency for
\[
u_t+F(D_x^2u,D_xu,u,t,x)=0.
\]
If $\halfup(0,\cdot)\le u_0\le\halfdown(0,\cdot)$ and comparison
holds for bounded sub- and supersolutions on every finite time
strip, then $\halfup=\halfdown=:u\in C(\overline Q)$ and
\[
\sup_{z\in K\cap\mathcal G_\rho}|u^\rho(z)-u(z)|
\longrightarrow0
\qquad\text{for every compact }K\subset\overline Q.
\]
\end{corollary}
\begin{proof}
The proof of \cref{thm:interior}, applied locally in space--time,
gives the sub- and supersolution properties in $Q$.
The initial inequalities and comparison on
$[0,T]\times\mathbb R^d$ for every $T>0$ yield
$\halfup=\halfdown=:u$, with $u(0,\cdot)=u_0$.
Continuity follows from the semicontinuity of the half-relaxed
limits, and local nodal convergence follows from the compactness
argument of \cref{cor:convergence}.
\end{proof}

\section{Verification of the local comparison hypothesis}
\label{sec:verification}

We give two complementary criteria for verifying \textnormal{(H3)}.
The first replaces the original nonmonotone operator with a nearby monotone auxiliary one, while keeping the numerical family suitably close. The second establishes a patchwise comparison property
for the original scheme. Throughout this section, \textnormal{(BS)}
denotes precisely \eqref{eq:BS}: monotonicity in the grid-function
argument with the scalar argument fixed.
\subsection{The transfer mechanism}

Assume (H1)--(H2). Transfer replaces the original numerical family
by a nearby one that approximately solves a monotone auxiliary scheme.
The latter admits standard localizing grid patches by locality and
grid density, without alignment to the original stencil.

\begin{definition}[Local averaging family]

A family of linear operators $J_\rho$ on grid functions is called a
\emph{local averaging family} if
\begin{equation}
 (J_\rho V)(x)=\sum_{y\in\mathcal G_\rho}j_\rho(x,y)V(y),
 \qquad j_\rho(x,y)\ge0,
 \qquad \sum_{y\in\mathcal G_\rho}j_\rho(x,y)=1,
\end{equation}
and, for every $K\Subset\Omega$,
\begin{equation}
 r_{\rho,K}:=
 \sup_{x\in K\cap\mathcal G_\rho}
 \sup_{\{y:\,j_\rho(x,y)>0\}}|y-x|_\star
 \longrightarrow0.
\end{equation}
The supremum is interpreted as zero when $K\cap\mathcal G_\rho$ is
empty. The identity $J_\rho=I$ is allowed.
\end{definition}

Local averaging preserves constants and is a contraction in the uniform
norm. Shrinking supports alone, however, do not imply
$J_\rho u^\rho-u^\rho\to0$ for an arbitrary bounded numerical family.

\noindent\textbf{(MTM) Monotone Transfer Mechanism.}

Let $J_\rho$ be a local averaging family and let $\mathcal M_\rho$ be a
local auxiliary operator satisfying \textnormal{(BS)}, with $M$ in
place of $S$ in \eqref{eq:BS}. Set $v^\rho=J_\rho u^\rho$.
These operators realize the \emph{monotone transfer mechanism},
abbreviated \textnormal{(MTM)}, if, for every $K\Subset\Omega$, the
following conditions hold:
\begin{enumerate}[label=\textnormal{(TM\arabic*)}]
\item \emph{Smooth-test compatibility.} For every
$\phi\in C_b^2(\mathbb R^d)$,
\begin{equation}
 \lim_{\substack{\rho\to0\\ \xi\to0}}
 \|\mathcal M_\rho[\phi+\xi]
     -\mathcal S_\rho[\phi+\xi]\|_K=0.
\end{equation}

\item \emph{Approximation of the numerical family.}
\begin{equation}
 \|v^\rho-u^\rho\|_K\longrightarrow0.
\end{equation}

\item \emph{Vanishing transfer defect.}
\begin{equation}
 \|\mathcal M_\rho[v^\rho]
     -J_\rho(\mathcal S_\rho[u^\rho])\|_K\longrightarrow0.
\end{equation}
\end{enumerate}

\begin{proposition}
\label[proposition]{prop:local-transfer}
If \textnormal{(MTM)} holds, then \textnormal{(H3)} holds with the
auxiliary pair $(\mathcal M_\rho,J_\rho u^\rho)$. In particular, this
pair satisfies \textnormal{(DC)}.
\end{proposition}

\begin{proof}
Condition \textnormal{(TM1)} transfers the relaxed consistency in
\textnormal{(H2)} to $\mathcal M_\rho$, while \textnormal{(TM2)} gives
the approximation required in \textnormal{(H3)}.
For $K\Subset\Omega$, locality of $J_\rho$ ensures that, for small
$\rho$, all averaging nodes used at points of $K$ lie in $\Omega$.
Thus \textnormal{(H1)} gives
$J_\rho(\mathcal S_\rho[u^\rho])=0$ on $K\cap\mathcal G_\rho$,
and \textnormal{(TM3)} implies
$\|\mathcal M_\rho[v^\rho]\|_K\to0$.
The residual may be extended arbitrarily to boundary nodes, since
these extensions do not enter this argument for sufficiently small
$\rho$.

Locality of $\mathcal M_\rho$ and grid density provide a localizing
family of grid balls: if $0<r<s<R$ and
$\overline B_R(x_0)\Subset\Omega$, then
$D=\mathcal G_\rho\cap B_s(x_0)$ contains the inner grid ball and its
stencil boundary lies in $B_R(x_0)\setminus B_s(x_0)$ for small $\rho$.
In an upper crossing, choose $z\in D$ maximizing
$v^\rho-\phi-c$, and let $m=v^\rho(z)-\phi(z)-c>0$.
Then $v^\rho\le\phi+c+m$ on the stencil of $z$, with equality at its
central node. Locality and \textnormal{(BS)} give
\[
 \mathcal M_\rho[v^\rho](z)
 \ge\mathcal M_\rho[\phi+c+m](z).
\]
In a lower crossing, maximize $\phi+c-v^\rho$ instead and obtain
\[
 \mathcal M_\rho[\phi+c-m](z)
 \ge\mathcal M_\rho[v^\rho](z).
\]
These are the two inequalities in \cref{def:dc}, with $\theta=m$.
\end{proof}

\begin{remark}[The identity transfer]
If $\mathcal S_\rho$ is monotone, taking $J_\rho=I$ and
$\mathcal M_\rho=\mathcal S_\rho$ gives (MTM) directly.
The general averaging formulation also permits a regularized family
$J_\rho u^\rho$ when smoother grid data are needed. 
\end{remark}

\begin{remark}[Transfer to a nonmonotone auxiliary operator]

Monotonicity of $\mathcal M_\rho$ is needed only to obtain (DC).
If (DC) holds directly for $(\mathcal M_\rho,J_\rho u^\rho)$,
(TM1)--(TM3) still give the more general transfer mechanism (TM):
$\textnormal{(MTM)}\Longrightarrow\textnormal{(TM)}
\Longrightarrow\textnormal{(H3)}$.
\end{remark}

\subsection{Patchwise comparison}

The direct route requires strict patchwise comparison (SPC) for
the original operator. It implies (DC), and hence (H3) under
(H1)--(H2).

\noindent\textbf{(SPC) Strict Patchwise Comparison.}

We say that $\mathcal S_\rho$ satisfies \textnormal{(SPC)}
if it admits a localizing family of patches $\mathfrak D^{\mathcal S}$
such that, for every $\rho$, every $D\in\mathfrak D^{\mathcal S}_\rho$,
every pair of grid functions $v,w$, and every $\delta>0$,
\[
    v\le w \quad\text{on }\partial_\rho D,
    \qquad
    \mathcal S_\rho[v]\le\mathcal S_\rho[w]-\delta
    \quad\text{on }D
\]
imply $v\le w$ on $D$.

\begin{proposition}
\label[proposition]{prop:all-routes}
Assume \textnormal{(H1)--(H2)}. Then
\[
 \textnormal{(SPC)}\Longrightarrow\textnormal{(H3)}.
\]
Under \textnormal{(SPC)}, the pair
$(\mathcal S_\rho,u^\rho)$ satisfies \textnormal{(DC)} with
$\theta=0$.
\end{proposition}

\begin{proof}
Suppose that \textnormal{(SPC)} holds and that $v\le w$ on
$\partial_\rho D$ but $\max_D(v-w)>0$. If
$\mathcal S_\rho[v](x)<\mathcal S_\rho[w](x)$ for every $x\in D$,
finiteness of $D$ gives
\[
 \delta:=\min_{x\in D}
 \bigl(\mathcal S_\rho[w](x)-\mathcal S_\rho[v](x)\bigr)>0.
\]
Property \textnormal{(SPC)} would imply $v\le w$ on $D$, a
contradiction. Thus there is $z\in D$ with
$\mathcal S_\rho[v](z)\ge\mathcal S_\rho[w](z)$.

Apply this observation to $(v,w)=(u^\rho,\phi+c)$ in the
upper-crossing configuration to obtain \eqref{eq:dc-upper}, and
to $(v,w)=(\phi+c,u^\rho)$ in the lower-crossing configuration
to obtain \eqref{eq:dc-lower}, both with $\theta=0$.
Hence $(\mathcal S_\rho,u^\rho)$ satisfies \textnormal{(DC)}. Taking
$\mathcal M_\rho=\mathcal S_\rho$ and $v^\rho=u^\rho$,
consistency follows from \textnormal{(H2)}, approximation is exact,
and the auxiliary equation is exact by \textnormal{(H1)}.
This proves \textnormal{(H3)}.
\end{proof}

\begin{proposition}[Comparison under constant shifts]

Suppose that
\begin{equation}
 \mathcal S_\rho[V+t](x)\ge\mathcal S_\rho[V](x)
 \qquad
 \text{for every grid function }V,\ t\ge0,
 \text{ and }x\in\mathcal G_\rho\cap\Omega.
 \label{eq:shift-monotonicity}
\end{equation}
Then
\[
 \textnormal{(BS)}\Longrightarrow\textnormal{(SPC)}.
\]
\end{proposition}

\begin{proof}
Suppose that \textnormal{(BS)} holds. Use the localizing family
of grid balls constructed in the proof of \cref{prop:local-transfer},
with the stencils of $\mathcal S_\rho$. Suppose that the hypotheses
of \textnormal{(SPC)} hold but $m=\max_D(v-w)>0$, and choose
$z\in D$ maximizing $v-w$. Then $v\le w+m$ on the stencil of $z$,
with equality at the central node. Locality, \textnormal{(BS)}, and
\eqref{eq:shift-monotonicity} give
\[
 \mathcal S_\rho[v](z)
 \ge\mathcal S_\rho[w+m](z)
 \ge\mathcal S_\rho[w](z),
\]
contradicting the strict operator gap.
\end{proof}

\begin{remark}
	
	A useful way to verify \textnormal{(SPC)} for Bellman and Isaacs schemes is
	through inverse positivity of the discrete operators obtained by freezing
	the controls. Related matrix-based discrete comparison arguments for
	nonmonotone elliptic discretizations are studied in
	\cite{PollockZhu2020Matrix}. Suppose that, on a patch \(D\), every row policy \(\pi\)
	defines an affine operator of the form
	\[
	\mathcal F_{\rho,D}^{\pi}[V]
	=
	A_{D,\rho}^{\pi}V_D
	+
	C_{D,\rho}^{\pi}V_{\partial_\rho D}
	-
	f_{D,\rho}^{\pi}.
	\]
	If, for every patch in a local patch basis and every admissible row policy,
	\begin{equation}
		(A_{D,\rho}^{\pi})^{-1}\ge0,
		\qquad
		-(A_{D,\rho}^{\pi})^{-1}C_{D,\rho}^{\pi}\ge0
		\label{eq:matrix-patch-comparison}
	\end{equation}
	entrywise, then the corresponding frozen problem is order preserving with
	respect to both the interior forcing and the stencil-boundary data.
	
	For Bellman and Isaacs operators, suitable rowwise choices of the controls
	reduce the comparison of two nonlinear residuals to such a frozen linear
	problem. Hence \eqref{eq:matrix-patch-comparison} yields
	\textnormal{(SPC)}, and therefore the discrete comparison condition required
	in \textnormal{(H3)}. The main scheme-dependent issue is thus to construct
	an appropriate local patch basis and to prove
	\eqref{eq:matrix-patch-comparison} uniformly over all row policies.
	This strategy is carried out explicitly for the one-dimensional \(Q^2\) Isaacs scheme below.
\end{remark}
\section{\texorpdfstring{Application of the transfer mechanism: CWENO for one-dimensional Hamilton--Jacobi equations}{Application of the transfer mechanism: CWENO for one-dimensional Hamilton-Jacobi equations}}
\label{sec:cweno}

For the one-dimensional semi-Lagrangian CWENO scheme of \cite{CarliniFerrettiPredaSemplice2025}, we establish convergence under the weaker mesh condition \(h=o(\tau)\), instead of \(h=O(\tau^2)\), and for bounded uniformly continuous initial data, without requiring Lipschitz continuity. The proof relies on comparison of the nonmonotone CWENO update with a monotone scheme.

Consider
\begin{equation}
\begin{cases}
u_t+H(u_x)=0,
    &(t,x)\in(0,\infty)\times\mathbb R,\\
u(0,x)=u_0(x),
    &x\in\mathbb R,
\end{cases}
\label{eq:cweno-pde}
\end{equation}
where $u_0\in\operatorname{BUC}(\mathbb R)$
(the space of bounded uniformly continuous functions),
$H\in C^2(\mathbb R)$, and
$H''\geq m_H>0$. Let $L=H^*$, where ${}^*$ denotes the Legendre--Fenchel transform, $x_i=ih$, and $t_n=n\tau$. 

In the autonomous setting considered in
\cite[Section~4 and Appendix~A]{CarliniFerrettiPredaSemplice2025}, the
scheme can be written as
\begin{equation}
U_i^{n+1}
=
\min_{a\in\mathbb R}
\left\{
R_h[U^n](x_i-\tau a)+\tau L(a)
\right\},
\qquad
U_i^0=u_0(x_i),
\label{eq:cweno-update}
\end{equation}
where, on $[x_j,x_{j+1}]$,
\[
R_h[V]=\omega_0P_0+\omega_LP_L+\omega_RP_R,
\qquad P_0=\frac43Q-\frac16(P_L+P_R).
\]
Here $Q,P_L,P_R$ are the Lagrange interpolants of $V$ on
$\{x_{j-1},x_j,x_{j+1},x_{j+2}\}$,
$\{x_{j-1},x_j,x_{j+1}\}$ and $\{x_j,x_{j+1},x_{j+2}\}$,
respectively. For $k\in\{0,L,R\}$,
\[
\omega_k=\frac{\alpha_k}{\alpha_0+\alpha_L+\alpha_R},
\qquad
\alpha_k=
\begin{cases}
d_k(h^2+\beta_k)^{-2}, & \text{CWENO},\\
d_k[1+(\Theta/(h^2+\beta_k))^2], & \text{CWENOZ},
\end{cases}
\]
with $(d_0,d_L,d_R)=(3/4,1/8,1/8)$,
$\Theta=|2\beta_0-\beta_L-\beta_R|$, and
\[
(\beta_0,\beta_L,\beta_R)
=(\mathcal I_j[Q],\mathcal I_j[P_L],\mathcal I_j[P_R]),
\qquad
\mathcal I_j[P]=\sum_{\ell=2}^{3}h^{2\ell-3}
\int_{x_j}^{x_{j+1}}\bigl(P^{(\ell)}(x)\bigr)^2\,dx;
\]
see \cite[Section~3.1]{CarliniFerrettiPredaSemplice2025}.

Let $T_R$ denote the update operator in \eqref{eq:cweno-update},
and define $T_I$ by replacing $R_h$ with the piecewise-linear
interpolant $I_h$. Both updates attain their minima for bounded
grid functions, since the reconstructions are
continuous and bounded and $L$ is superlinear. Moreover,
$L\in C^2(\mathbb R)$ and $0<L''\le1/m_H$.
Set $\delta_j^2V=V_{j+1}-2V_j+V_{j-1}$, and let
$\omega_{u_0}(r)=\sup_{|x-y|\le r}|u_0(x)-u_0(y)|$.

The key comparison estimate is the following.
\begin{lemma}[One-sided reconstruction and comparison of updates]
\label[lemma]{lem:cweno-one-sided}
If $V$ is bounded and $\delta_j^2V\le Kh^2$ for all $j$, with $K\ge0$, then
\begin{equation}
R_h[V]\ge I_h[V]-\frac{Kh^2}{8},
\qquad
T_IV-\frac{Kh^2}{8}\le T_RV
\le T_IV+\frac{h^2}{8m_H\tau}.
\end{equation}
\end{lemma}

The proof is given in Appendix~\ref{app:cweno-proof}.

The reconstruction identity \eqref{eq:cweno-signed-reconstruction},
established there, also yields
\begin{equation}
\|R_h[\phi]-I_h[\phi]\|_\infty
\le\frac{h^2}{8}\|\phi_{xx}\|_\infty,
\qquad
\|R_h[U^0]-I_h[U^0]\|_\infty\le\frac{\omega_{u_0}(h)}4,
\label{eq:cweno-test-defect}
\end{equation}
and $\|R_h[V]\|_\infty\le\tfrac32\|V\|_\infty$ for bounded $V$.

To specify convergence away from the nodes, define
\[
\widetilde U_{h,\tau}(t,x)=R_h[U^n](x),
\qquad t_n\le t<t_{n+1}.
\]

\begin{theorem}[CWENO convergence under a relaxed mesh condition]
\label[theorem]{thm:cweno-convergence}
Under the preceding assumptions and $h=o(\tau)$, the scheme
\eqref{eq:cweno-update} realizes \textnormal{(MTM)} with $J_{h,\tau}=I$
on each finite time strip. The extensions
$\widetilde U_{h,\tau}$ converge locally
uniformly on $[0,T]\times\mathbb R$, for every $T>0$, to the bounded
uniformly continuous viscosity solution of \eqref{eq:cweno-pde}.
\end{theorem}

\begin{proof}
First, every update is semiconcave at the grid nodes:
\begin{equation}
\delta_j^2U^n\le\frac{h^2}{m_H\tau},\qquad n\ge1.
\label{eq:cweno-semiconcavity}
\end{equation}
To see this, use a minimizing foot $y_*$ for $(T_RV)_j$ as a competitor
at $x_j-h$ and $x_j+h$. The reconstruction values cancel in the second
difference, and $L''\le1/m_H$ proves \eqref{eq:cweno-semiconcavity}
for any bounded input $V$. Consequently, \cref{lem:cweno-one-sided}
gives
\begin{equation}
\|T_RU^n-T_IU^n\|_\infty\le\frac{h^2}{8m_H\tau},
\qquad n\ge1.
\label{eq:cweno-family-defect}
\end{equation}

For stability and the initial condition, let $W^0=U^0$ and
$W^{n+1}=T_IW^n$. Monotonicity and preservation of constant differences
make $T_I$ nonexpansive in the uniform norm, and comparison with
constants gives
$\|W^n\|_\infty\le\|u_0\|_\infty+t_n|H(0)|$.
The first-step error is at most $\omega_{u_0}(h)/4$ by
\eqref{eq:cweno-test-defect}. Iterating \eqref{eq:cweno-family-defect}
therefore gives
\begin{equation}
\max_{n\tau\le T}\|U^n-W^n\|_\infty
\le\frac{\omega_{u_0}(h)}4+\frac{Th^2}{8m_H\tau^2}
\longrightarrow0.
\label{eq:cweno-stability-transfer}
\end{equation}
This proves a uniform bound $\|U^n\|_\infty\le M_T$. The family $W$
inherits $u_0$ uniformly at $t=0$: approximate $u_0$ uniformly by a
bounded Lipschitz function $v_0$. If $P$ is its Lipschitz constant,
the bounds $I_h[v_0](x_i-\tau a)\ge v_0(x_i)-P\tau|a|$ and the
competitor $a=0$ give $\|T_Iv_0-v_0\|_\infty\le C_P\tau$.
Nonexpansivity then gives
$\|W^n-u_0\|_\infty\le2\|v_0-u_0\|_\infty+C_Pt_n$.
Together with \eqref{eq:cweno-stability-transfer}, this also proves
inheritance of the initial condition for $U$.

We next localize the controls. For either reconstruction applied to
$U^n$, comparison of a minimizer with $a=0$ gives
$\tau L(a)\le\tfrac52M_T+\tau L(0)$.
Choose $C_T>\tfrac52M_T$ and set
\[
A_\tau=\max\{|a|:L(a)\le L(0)+C_T/\tau\},
\qquad \mathcal A_\tau=[-A_\tau,A_\tau].
\]
Superlinearity implies $A_\tau\to\infty$ and $\tau A_\tau\to0$.
Restricting both updates to $\mathcal A_\tau$ therefore preserves their
values on $U^n$, while their stencil radius is at most
$\tau A_\tau+2h\to0$. Define the normalized operators by
\begin{align*}
\mathcal S^{\rm C}_{h,\tau}[V](t_{n+1},x_i)
&=\frac{V_i^{n+1}-\min_{a\in\mathcal A_\tau}
\{R_h[V^n](x_i-\tau a)+\tau L(a)\}}{\tau},\\
\mathcal M_{h,\tau}[V](t_{n+1},x_i)
&=\frac{V_i^{n+1}-\min_{a\in\mathcal A_\tau}
\{I_h[V^n](x_i-\tau a)+\tau L(a)\}}{\tau}.
\end{align*}
The latter is monotone in the Barles--Souganidis sense.
For a fixed bounded smooth test $\phi$, minimizing velocities remain
in a fixed compact set: by \eqref{eq:cweno-test-defect}, they satisfy
\[
L(a)-\|\phi_x\|_\infty|a|
\le L(0)+C h^2/\tau.
\]
Since $h^2/\tau\to0$, superlinearity proves this claim. Taylor expansion
on that fixed set gives consistency of both normalized operators with
$\phi_t+H(\phi_x)$, including vanishing constant shifts.

Take $J_{h,\tau}=I$ and $v^{h,\tau}=U$. For \textnormal{(TM1)},
\eqref{eq:cweno-test-defect} and the uniform-norm nonexpansivity of
the minimum give
\[
\|\mathcal M_{h,\tau}[\phi+\xi]
-\mathcal S^{\rm C}_{h,\tau}[\phi+\xi]\|_K
\le\frac{h^2}{8\tau}\|\phi_{xx}\|_\infty\longrightarrow0,
\]
uniformly in $\xi$. Condition \textnormal{(TM2)} is an identity. For
\textnormal{(TM3)}, \eqref{eq:cweno-family-defect} gives
\begin{equation}
\|\mathcal M_{h,\tau}[U]-\mathcal S^{\rm C}_{h,\tau}[U]\|_K
\le\frac{h^2}{8m_H\tau^2}\longrightarrow0.
\end{equation}
Every $K\Subset(0,T)\times\mathbb R$ excludes the first time step
for sufficiently small $\tau$; that step has already been treated in
\eqref{eq:cweno-stability-transfer}. Thus \textnormal{(H1)--(H2)} hold,
\cref{prop:local-transfer} gives \textnormal{(H3)}, and
\cref{cor:cauchy} yields nodal convergence. Applying the argument on
a slightly longer time strip includes the time $T$ in the conclusion.

For $x\in[x_j,x_{j+1}]$,
\eqref{eq:cweno-signed-reconstruction} gives
\[
|R_h[U^n](x)-I_h[U^n](x)|
\le\frac18\max\{|\delta_j^2U^n|,|\delta_{j+1}^2U^n|\}.
\]
Nodal convergence to the continuous function $u$ makes these second
differences tend to zero locally uniformly. The piecewise-linear
spatial interpolants, extended constantly on each time interval,
converge to $u$ by uniform continuity on compact sets. Hence
$\widetilde U_{h,\tau}\to u$ locally uniformly as claimed.
\end{proof}

\section{\texorpdfstring{Application of patchwise comparison: $Q^2$ for one-dimensional Isaacs equations}{Application of patchwise comparison: Q2 for one-dimensional Isaacs equations}}
\label{sec:q2}
We apply the direct comparison mechanism to a one-dimensional
$Q^2$ scheme. The key observation is that eliminating the odd-index
unknowns produces a tridiagonal $M$-matrix on the even-index nodes.
Throughout this application, use $|x|_\star=|x|$ for the norm in \eqref{eq:density}.

Let $\mathcal A,\mathcal B$ be finite nonempty control sets and consider
\begin{equation}
\begin{cases}
\displaystyle\max_{\alpha\in\mathcal A}\min_{\beta\in\mathcal B}
\{-a^{\alpha\beta}(x)u''+b^{\alpha\beta}(x)u'+u-f^{\alpha\beta}(x)\}=0,
&x\in(0,1),\\
u(0)=g_0,\qquad u(1)=g_1.
\end{cases}
\label{eq:q2-pde}
\end{equation}
Assume comparison for this continuous Dirichlet problem and
$a^{\alpha\beta},b^{\alpha\beta},f^{\alpha\beta}\in C([0,1])$, with
\begin{equation}
0<a_0\le a^{\alpha\beta}(x)\le a_1,
\qquad |b^{\alpha\beta}(x)|\le B.
\label{eq:q2-coefficients}
\end{equation}
On $x_i=ih$, $h=1/(2M)$, use the one-dimensional $Q^2$ stencil
from the finite-difference realization studied in \cite{LiZhang2020}:
\begin{equation}
(\Qh U)_i=
\begin{cases}
\displaystyle\frac{2U_i-U_{i-1}-U_{i+1}}{h^2},&i\ \text{odd},\\[2mm]
\displaystyle\frac{U_{i-2}-8U_{i-1}+14U_i-8U_{i+1}+U_{i+2}}{4h^2},
&i\ \text{even}.
\end{cases}
\label{eq:Q2}
\end{equation}
Odd nodes are element centres and even nodes are element endpoints;
the element width is $2h$. With
$(D_h^0U)_i=(U_{i+1}-U_{i-1})/(2h)$, the scheme is
\begin{equation}
\mathcal S_{h,i}[U]
=\max_\alpha\min_\beta
\{a_i^{\alpha\beta}(\Qh U)_i+b_i^{\alpha\beta}(D_h^0U)_i
+U_i-f_i^{\alpha\beta}\}=0,
\quad 1\le i\le2M-1,
\label{eq:q2-scheme}
\end{equation}
with $U_0=g_0$, $U_{2M}=g_1$.

\begin{remark}[Nonmonotonicity and the choice of patches]
The positive distance-$2h$ coefficient at each even node violates
\eqref{eq:BS}. It also prevents uniform near-monotonicity on bounded
grid functions: for interior even-index nodes $x_i$ tending to $1/2$,
put $V=0$ and $W=h^2\mathbf1_{\{x_{i+2}\}}$. Each affine expression increases
by at least $a_0/4$ at $x_i$, so
$\mathcal S_{h,i}[W]-\mathcal S_{h,i}[V]\ge a_0/4$.
If $\mathcal S_h=\mathcal T_h+\mathcal P_h$ with $\mathcal T_h$
satisfying \eqref{eq:BS}, the central values agree and $V\le W$ gives
\[
a_0/4\le|\mathcal P_h[W](x_i)|+|\mathcal P_h[V](x_i)|.
\]
This rules out a perturbation tending uniformly to zero on bounded
grid functions.

Even patchwise comparison fails on a single even node. For one control
with $a=1$ and $b=f=0$, take $D=\{x_i\}$, $V=0$, and prescribe
$W_{i-2}=1$, $W_i=-1/(28+8h^2)$, all other stencil values being zero.
Then $V\le W$ on $\partial_hD$ and
$\mathcal S_{h,i}[V]=0<\mathcal S_{h,i}[W]=1/(8h^2)$, yet $V_i>W_i$.
Including the two neighbouring odd nodes changes the boundary response.
On $D=\{x_{i-1},x_i,x_{i+1}\}$, the homogeneous equations give
\[
U_i=\frac{6-h^2}{12+22h^2+4h^4}(U_{i-2}+U_{i+2}),
\qquad
U_{i\pm1}=\frac{U_i+U_{i\pm2}}{2+h^2}.
\]
For sufficiently small $h$, all these boundary coefficients are
positive. Eliminating the odd nodes is therefore the natural way
to recover comparison.
\end{remark}

Accordingly, take the aligned intervals
\begin{equation}
D_{p,q}^h=\{x_{2p+1},\ldots,x_{2q-1}\},
\qquad \partial_hD_{p,q}^h=\{x_{2p},x_{2q}\},
\quad 0\le p<q\le M.
\end{equation}
Their endpoints have even indices, so the stencils do not cross the
endpoints and the intermediate odd nodes remain inside the patch.
To check \eqref{eq:basis}, fix $x_0$ and $0<r<R$ with
$[x_0-R,x_0+R]\Subset(0,1)$. Choose even-index endpoints
\[
x_{2p}\in(x_0-R,x_0-r),\qquad
x_{2q}\in(x_0+r,x_0+R).
\]
Both choices are possible when $2h<R-r$, independently of the position
of $x_0$ relative to the grid. The resulting patch contains every node
of $B_r(x_0)$, lies in $B_R(x_0)$, and its two stencil-boundary nodes
belong to $B_R(x_0)\setminus B_r(x_0)$. Thus these intervals form a
localizing family.

\begin{theorem}[Convergence of the one-dimensional $Q^2$ Isaacs scheme]
\label{thm:q2-convergence}
Under the preceding assumptions, if
\begin{equation}
h^2+6hB<6a_0,
\label{eq:q2-mesh-condition}
\end{equation}
the scheme \eqref{eq:q2-scheme} satisfies \textnormal{(SPC)} on the
aligned intervals and has a unique solution $U^h$ with the prescribed
endpoint values. As $h\to0$,
\[
\max_{0\le i\le2M}|U_i^h-u(x_i)|\longrightarrow0,
\]
where $u\in C([0,1])$ is the viscosity solution of \eqref{eq:q2-pde}.
\end{theorem}

A row policy $\pi$ chooses a pair $(\alpha_i,\beta_i)$ at each node
of $D=D_{p,q}^h$. Write $A_D^\pi$ for its interior coefficient matrix
and $C_D^\pi$ for its two boundary columns.
\begin{lemma}[Policywise inverse positivity]
\label[lemma]{lem:q2-policy-positivity}
Under \eqref{eq:q2-mesh-condition}, every aligned interval and every
row policy satisfy
\begin{equation}
(A_D^\pi)^{-1}\ge0,\qquad
-(A_D^\pi)^{-1}C_D^\pi\ge0.
\label{eq:q2-policy-positivity}
\end{equation}
\end{lemma}
The proof in Appendix~\ref{app:q2-proof} eliminates the odd unknowns.
The resulting even-node coefficients are nonpositive off the diagonal
and strictly dominant on the diagonal.

\begin{proof}[Proof of \Cref{thm:q2-convergence}]
Write $L_{h,i}^{\alpha\beta}$ for the affine operator associated with
the controls $(\alpha,\beta)$ in
\eqref{eq:q2-scheme}. At each interior node choose
\[
\alpha_W\in\operatorname*{arg\,max}_\alpha\min_\beta
L_{h,i}^{\alpha\beta}[W],
\qquad
\beta_V\in\operatorname*{arg\,min}_\beta
L_{h,i}^{\alpha_W\beta}[V].
\]
Then
\[
\begin{aligned}
\mathcal S_{h,i}[V]&\ge\min_\beta L_{h,i}^{\alpha_W\beta}[V]
=L_{h,i}^{\alpha_W\beta_V}[V],\\
\mathcal S_{h,i}[W]&=\min_\beta L_{h,i}^{\alpha_W\beta}[W]
\le L_{h,i}^{\alpha_W\beta_V}[W].
\end{aligned}
\]
The asymmetric choices give the directions needed on subtraction.
The forcing cancels, and the resulting row policy satisfies
\[
A_D^\pi(V-W)_D+C_D^\pi(V-W)_{\partial_hD}
\le\mathcal S_h[V]-\mathcal S_h[W]\le0.
\]
Multiplication by the nonnegative interior inverse gives
\[
(V-W)_D\le
-(A_D^\pi)^{-1}C_D^\pi(V-W)_{\partial_hD}\le0.
\]
The last inequality uses the nonnegative boundary-transfer matrix
in \eqref{eq:q2-policy-positivity} and the ordering on $\partial_hD$.
We thus obtained \textnormal{(SPC)} on the localizing
family and uniqueness on the full interval $D=D_{0,M}^h$.
We prove existence by two finite Bellman policy iterations. 
For a fixed rowwise outer policy $\boldsymbol\alpha=(\alpha_i)_{x_i\in D}$,
define the inner Bellman operator
\[
F_h^{\boldsymbol\alpha}[V]_i
=\min_{\beta\in\mathcal B}L_{h,i}^{\alpha_i\beta}[V].
\]

First solve $F_h^{\boldsymbol\alpha}[Z]=0$ with the prescribed endpoints.
Start with the solution $Z^{(0)}$ of any inner row policy. If $Z^{(k)}$
solves a policy equation, then $F_h^{\boldsymbol\alpha}[Z^{(k)}]\le0$.
Choose $\boldsymbol\beta^{(k+1)}$ attaining the minimum at $Z^{(k)}$
and solve its linear equations for $Z^{(k+1)}$. Their difference gives
\[
A_D^{\boldsymbol\alpha,\boldsymbol\beta^{(k+1)}}
(Z^{(k+1)}-Z^{(k)})_D
=-F_h^{\boldsymbol\alpha}[Z^{(k)}]\ge0.
\]
There is no boundary contribution because both iterates have the same
endpoint values. Inverse positivity implies $Z^{(k+1)}\ge Z^{(k)}$.
Only finitely many inner row policies exist, each with a unique linear
solution; a monotone sequence in this finite set must become stationary.
The displayed identity then gives $F_h^{\boldsymbol\alpha}[Z]=0$.
Comparison gives uniqueness of this inner Bellman solution.

Now start with any outer row policy $\boldsymbol\alpha^{(0)}$ and let
$U^{(m)}$ be the unique inner Bellman solution for
$\boldsymbol\alpha^{(m)}$. Since the full operator maximizes over the
outer control, $\mathcal S_h[U^{(m)}]\ge0$. Select
$\boldsymbol\alpha^{(m+1)}$ attaining this maximum at every node and
solve the new inner Bellman problem completely. Then
\[
F_h^{\boldsymbol\alpha^{(m+1)}}[U^{(m)}]
=\mathcal S_h[U^{(m)}]\ge0,
\qquad
F_h^{\boldsymbol\alpha^{(m+1)}}[U^{(m+1)}]=0.
\]
Comparison for the fixed-outer-policy Bellman operator implies
$U^{(m+1)}\le U^{(m)}$. There are only finitely many outer row policies
and one inner solution for each, so this sequence also becomes stationary.
At a stationary step the two identities give $\mathcal S_h[U]=0$.
Thus the nested procedure proves existence in finitely many policy solves.
The monotonicity argument uses a complete inner solve before each outer
update; it does not require simultaneous improvement of both controls.

Put $F_0=\max_{\alpha,\beta}\|f^{\alpha\beta}\|_\infty$ and
$g_* =\max\{|g_0|,|g_1|\}$. Stability follows by comparison with
constants. For $C_0=\max\{F_0,g_*\}$, every affine operator satisfies
\[
L_{h,i}^{\alpha\beta}[C_0]=C_0-f_i^{\alpha\beta}\ge0,
\qquad
L_{h,i}^{\alpha\beta}[-C_0]=-C_0-f_i^{\alpha\beta}\le0.
\]
The endpoint values lie between $-C_0$ and $C_0$, so comparison gives
$-C_0\le U^h\le C_0$ and hence $\|U^h\|_\infty\le C_0$.\par
 We now construct upper and lower barriers to show that the numerical solutions approach the prescribed Dirichlet values. We start with the affine function \(\ell(x)=(1-x)g_0+xg_1\), which agrees with the boundary data, and consider barriers of the form \(\ell\pm C\psi\). Here \(\psi\) is nonnegative, vanishes at both endpoints, and is chosen concave so that the diffusion dominates the drift, i.e.
$a_0(-\psi'')>B|\psi'|$.
 The second derivative of an exponential
gains a factor $\kappa$ relative to its first derivative, so symmetrize
about $1/2$ and subtract the resulting
profile from its endpoint value:
\[
\psi(x)=\cosh(\kappa/2)-\cosh\bigl(\kappa(x-1/2)\bigr),
\qquad \kappa>B/a_0.
\]
Writing $\theta=\kappa(x-1/2)$ and $\mu=\kappa(a_0\kappa-B)>0$,
every pair of controls satisfies
\[
-a^{\alpha\beta}\psi''+b^{\alpha\beta}\psi'
\ge a_0\kappa^2\cosh\theta-B\kappa|\sinh\theta|
\ge\mu\cosh\theta\ge\mu.
\]
Here $|\sinh\theta|\le\cosh\theta$ explains the drift domination.
The parity of the stencil causes no difficulty: for a smooth profile
$g$, the symmetric Taylor expansions give
\[
(\Qh g)_i=
\begin{cases}
-g''(x_i)-\dfrac{h^2}{12}g^{(4)}(x_i)+O(h^4),&i\ \text{odd},\\[1mm]
-g''(x_i)+\dfrac{h^2}{6}g^{(4)}(x_i)+O(h^4),&i\ \text{even},
\end{cases}
\qquad D_h^0g=g'+O(h^2).
\]
Indeed, $g(x_i-h)+g(x_i+h)=2g+h^2g''+h^4g^{(4)}/12+O(h^6)$;
the analogous sum at distance $2h$ has coefficients $4$ and $4/3$
in place of $1$ and $1/12$. Inserting these sums in \eqref{eq:Q2}
gives the two displayed formulas. For the smooth barrier the
remainders are uniform, including nodes next to the endpoints.
Thus, for all sufficiently small $h$,
\[
a_i^{\alpha\beta}(\Qh\psi)_i+b_i^{\alpha\beta}(D_h^0\psi)_i
+\psi(x_i)\ge\mu/2
\]
uniformly in the controls. Since $\Qh\ell=0$, $D_h^0\ell=g_1-g_0$
and $\|\ell\|_\infty\le g_*$, choose
\[
K=\frac{2(F_0+B|g_1-g_0|+g_*)}{\mu}.
\]
For every pair of controls $(\alpha,\beta)$, the corresponding affine discrete operator  is nonnegative
on $\ell+K\psi$ and nonpositive on $\ell-K\psi$. Hence,  they are a supersolution
and a subsolution with the prescribed endpoint values. Comparison gives
\[
|U_i^h-\ell(x_i)|\le K\psi(x_i)
\le K\kappa\sinh(\kappa/2)\min\{x_i,1-x_i\}.
\]
This uniform endpoint decay ensures the Dirichlet condition for both
half-relaxed limits.

Finally, the stencil radius is $2h$, and second-order Taylor expansion
gives $\Qh\phi\to-\phi''$ and $D_h^0\phi\to\phi'$ locally uniformly
for every $C^2$ test function. The finite control sets and continuous
coefficients yield consistency of the max--min operator. Constant
shifts satisfy $\mathcal S_h[\phi+\xi]=\mathcal S_h[\phi]+\xi$.
This proves \textnormal{(H1)--(H2)}, and \textnormal{(SPC)} gives
\textnormal{(H3)} by \cref{prop:all-routes}. The convergence criterion
\cref{thm:main} completes the proof.
\end{proof}

\section{Numerical simulations}
\label{sec:numerics}

We test the CWENO and one-dimensional $Q^2$ Isaacs schemes of
\cref{sec:cweno,sec:q2}, respectively, focusing on the transfer
estimate for non-Lipschitz data and on max--min control switching and mesh refinement.

\subsection{\texorpdfstring{A Hamilton--Jacobi problem with non-Lipschitz initial data}{A Hamilton-Jacobi problem with non-Lipschitz initial data}}

We consider the one-periodic version of the Hamilton--Jacobi Cauchy
problem \eqref{eq:cweno-pde}, with $H(p)=p^2/2$, final time $T=1/8$,
and initial datum
\[
u_0(x)=\sqrt{\operatorname{dist}(x-1/3,\mathbb Z)}.
\]
This datum is bounded and uniformly continuous, but not Lipschitz.
For $k\in\{4,8,16,20\}$, we use a periodic grid with
\[
N=k^3,\qquad h=N^{-1}=k^{-3},\qquad \tau=k^{-2}.
\]
Thus $h=\tau^{3/2}=o(\tau)$, as required by
\cref{thm:cweno-convergence}. We use the nonlinear weights of
\cref{sec:cweno} and evaluate the reference solution by the periodic Hopf--Lax formula.

\begin{figure}[htbp]
    \centering
    \includegraphics[width=\linewidth]{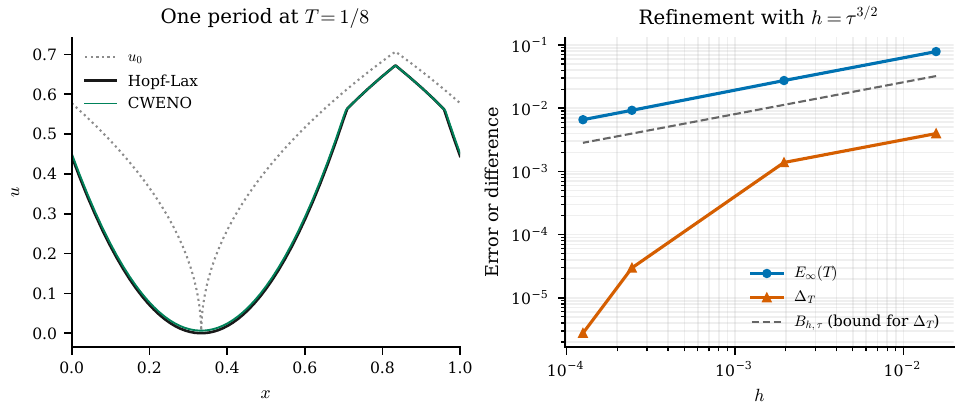}
    \caption{CWENO experiment for non-Lipschitz initial data with
    $h=\tau^{3/2}$. Left: initial datum and exact and CWENO profiles
    at $T=1/8$ on the finest grid ($k=20$). Right: final-time maximum
    error, difference from the monotone companion, and the bound
    from \eqref{eq:cweno-stability-transfer}.}
    \label{fig:cweno-numerical}
\end{figure}

Let $U^n$ denote the CWENO solution and $W^n$ the solution of its
piecewise-linear, monotone companion introduced in the proof of
\cref{thm:cweno-convergence}. We measure the final-time maximum error
and the difference between the two numerical families by
\[
E_\infty(T)=\max_i|U_i^{T/\tau}-u(T,x_i)|,
\qquad
\Delta_T=\max_{0\le n\tau\le T}\|U^n-W^n\|_\infty.
\]
For this datum, $\omega_{u_0}(h)\le\sqrt h$ and $m_H=1$.
The estimate \eqref{eq:cweno-stability-transfer} therefore gives
\[
\Delta_T\le B_{h,\tau}
:=\frac{\sqrt h}{4}+\frac{Th^2}{8\tau^2}.
\]
\Cref{fig:cweno-numerical} compares the initial, exact, and CWENO
profiles at the final time, and displays $E_\infty(T)$, $\Delta_T$,
and $B_{h,\tau}$ under refinement. Both measured quantities decrease,
and the difference from the monotone companion remains below the
predicted bound. This illustrates the transfer estimate used in the
convergence proof without a uniform numerical Lipschitz assumption.

\subsection{An Isaacs problem with switching controls}

We test \eqref{eq:q2-scheme} with two controls for each player,
$r,s\in\{-1,1\}$, variable diffusion, nonzero drift and nonzero
Dirichlet data. Choose
\[
a^{r,s}(x)=1+\frac15\sin(2\pi x)+\frac r5+\frac s{10},
\qquad b^{r,s}=\frac25r-\frac3{10}s+\frac35rs.
\]
For the coefficients in \eqref{eq:q2-coefficients}, we have $a_0=1/2$, $a_1=3/2$, and $B=13/10$; the sufficient mesh
condition in \eqref{eq:q2-mesh-condition} is $h^2+7.8h<3$. In particular, the diffusion contrast is
$3$, and the one-dimensional positivity proof applies without any
additional contrast bound.

Take the smooth exact solution
\[
u_{\rm ex}(x)=0.3-0.5x+\sin(2\pi x),
\qquad g_0=0.3,\quad g_1=-0.2,
\]
and a common forcing computed from its continuous derivatives:
\[
f(x)=\max_{r\in\{-1,1\}}\min_{s\in\{-1,1\}}
\{-a^{r,s}(x)u_{\rm ex}''(x)+b^{r,s}u_{\rm ex}'(x)+u_{\rm ex}(x)\}.
\]
The forcing is continuous and independent of the grid. Its derivative
may change where controls switch; the exact solution remains smooth.
We use $h=1/n$ for the six even values of $n$ in
\cref{tab:q2-numerical}. All satisfy the mesh condition.

The equations are solved by nested policy iteration: for each outer
policy, stabilize the inner minimizing policy, then update the outer
maximizer, retaining controls at numerical ties. Each linear policy
problem is solved by the odd-node elimination of
Appendix~\ref{app:q2-proof} and a tridiagonal solve. Stopping is checked
against the full max--min residual in the original stencils, with
tolerance $2\times10^{-9}$. As an additional check, we assemble the discrete operators corresponding to the four pairs of controls and evaluate the nonlinear Isaacs operator at the computed solution. We then calculate the maximum absolute value of this operator over the interior grid nodes. The largest value obtained over all the grids is \QTwoMaximumResidual.

\begin{table}[htbp]

\centering
\begin{tabular}{rrrr}
\hline
$n$ & $E_\infty$ & Order & $R_\infty$\\
\hline
16 & $4.73\times10^{-3}$ & -- & $2.34\times10^{-13}$ \\
32 & $1.28\times10^{-3}$ & 1.89 & $1.10\times10^{-12}$ \\
64 & $2.94\times10^{-4}$ & 2.12 & $4.06\times10^{-12}$ \\
128 & $7.41\times10^{-5}$ & 1.99 & $1.62\times10^{-11}$ \\
256 & $1.88\times10^{-5}$ & 1.98 & $8.53\times10^{-11}$ \\
512 & $4.76\times10^{-6}$ & 1.98 & $4.08\times10^{-10}$ \\
\hline
\end{tabular}

\caption{One-dimensional $Q^2$ Isaacs test. Here
$E_\infty=\max_i|U_i^h-u_{\rm ex}(x_i)|$,
the observed order is $\log_2(E_\infty(2h)/E_\infty(h))$, and
$R_\infty=\max_{1\le i\le n-1}|\mathcal S_{h,i}[U^h]|$ is
evaluated in the original stencils.}
\label{tab:q2-numerical}
\end{table}

\begin{figure}[htbp]

\centering
\includegraphics[width=\linewidth]{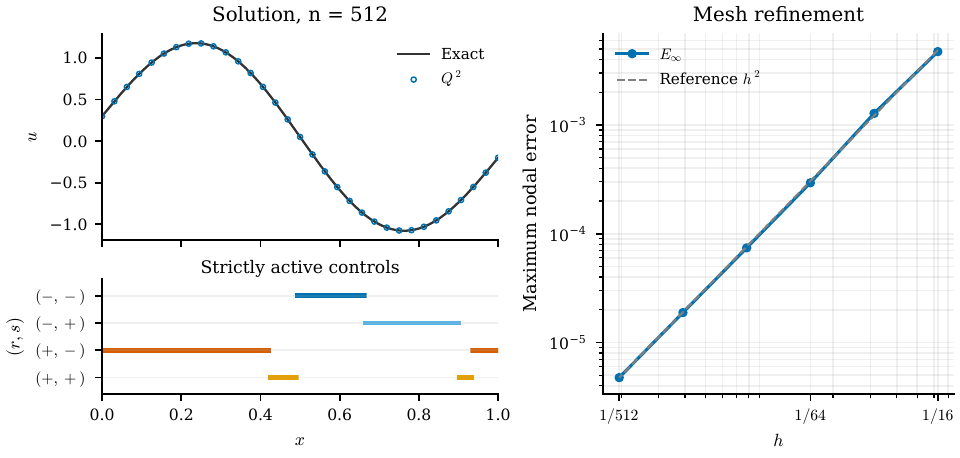}
\caption{One-dimensional Isaacs experiment. Left: exact solution and
sampled numerical values on $n=\QTwoFinestN$, with the strictly
active pairs $(r,s)$ below. Right: maximum nodal error under refinement
and a reference line proportional to $h^2$.}
\label{fig:q2-numerical}
\end{figure}

The maximum error decreases from \QTwoFirstError\ to
\QTwoLastError, with final observed order \QTwoLastOrder.
These results are consistent with second-order convergence for this
test; they do not assert an accuracy rate for general continuous data.
\Cref{fig:q2-numerical} shows the solution, the selected controls and
the refinement errors.

The game involves both optimizations. Declare a selected pair strictly
active only when the outer maximizing gap and its inner minimizing
gap both exceed $10^{-6}$. On $n=\QTwoFinestN$, the pairs
$(-1,-1),(-1,1),(1,-1),(1,1)$ are strictly active on
\QTwoStrictFractions\ of the interior nodes, respectively.
Moreover, evaluating $\min_s\max_r L_h^{r,s}[U^h]$ instead of
$\max_r\min_s L_h^{r,s}[U^h]$ gives a positive difference greater
than $10^{-6}$ on \QTwoOrderGapFraction\ of these nodes, with maximum
\QTwoOrderGap. Thus the order of the two optimizations matters in this
test. 

\appendix
\section{Technical proofs}
\label{app:technical-proofs}

\subsection{CWENO reconstruction and update comparison}
\label{app:cweno-proof}

\begin{proof}[Proof of \Cref{lem:cweno-one-sided}]
For $x=x_j+sh$, $0\le s\le1$, we have
$I_h[V](x)=(1-s)V_j+sV_{j+1}$.
The Newton interpolation formulas for the polynomials defined in
\cref{sec:cweno} give
\begin{align*}
P_L(x)-I_h[V](x)
&=\frac{s(s-1)}2\,\delta_j^2V,\\
P_R(x)-I_h[V](x)
&=\frac{s(s-1)}2\,\delta_{j+1}^2V,\\
Q(x)-P_L(x)
&=\frac{(s+1)s(s-1)}6
  \bigl(\delta_{j+1}^2V-\delta_j^2V\bigr).
\end{align*}
Combining the first and third identities yields
\[
Q(x)-I_h[V](x)
=-\frac{s(1-s)}2
\left(
\frac{2-s}{3}\,\delta_j^2V
+\frac{1+s}{3}\,\delta_{j+1}^2V
\right).
\]
Using the definition of $P_0$, we obtain
\[
P_0(x)-I_h[V](x)
=-\frac{s(1-s)}2
\left(
\frac{13-8s}{18}\,\delta_j^2V
+\frac{5+8s}{18}\,\delta_{j+1}^2V
\right).
\]
For both CWENO and CWENOZ, the nonlinear weights are nonnegative
and sum to one. The definition of $R_h$ and the preceding identities
therefore give
\begin{equation}
R_h[V](x)-I_h[V](x)
=-\frac{s(1-s)}2
\bigl(
\lambda_j(s)\delta_j^2V
+(1-\lambda_j(s))\delta_{j+1}^2V
\bigr),
\label{eq:cweno-signed-reconstruction}
\end{equation}
where
\[
\lambda_j(s)=\omega_L+\omega_0\frac{13-8s}{18},
\qquad
1-\lambda_j(s)=\omega_R+\omega_0\frac{5+8s}{18}.
\]
Both coefficients are nonnegative for $0\le s\le1$ and sum to one.
The semiconcavity assumption therefore yields
\[
R_h[V](x)-I_h[V](x)
\ge-\frac{s(1-s)}2\,Kh^2
\ge-\frac{Kh^2}{8}.
\]
Adding $\tau L((x_i-y)/\tau)$ to this inequality at $x=y$
and minimizing over $y\in\mathbb R$ gives
$(T_RV)_i\ge(T_IV)_i-Kh^2/8$.

For the upper bound, let $y_*$ minimize
$F_i(y)=I_h[V](y)+\tau L((x_i-y)/\tau)$. If $y_*$ is a node,
interpolation gives $(T_RV)_i\le(T_IV)_i$. Otherwise $y_*$ lies
inside a cell, where $F_i'(y_*)=0$ and
$F_i''\le1/(m_H\tau)$. At the closest endpoint $x_k$, whose distance
from $y_*$ is at most $h/2$, Taylor's formula gives
\[
(T_RV)_i\le V_k+\tau L((x_i-x_k)/\tau)
=F_i(x_k)\le F_i(y_*)+\frac{h^2}{8m_H\tau}.
\]
\end{proof}

\subsection{Odd-node elimination for the one-dimensional policy matrices}
\label{app:q2-proof}
\begin{proof}[Proof of \Cref{lem:q2-policy-positivity}]
Fix an aligned interval $D_{p,q}^h$ and any row policy. Suppress its
superscript and write the frozen equations with an arbitrary forcing:
\[
a_i(\Qh U)_i+b_i(D_h^0U)_i+U_i=d_i,
\qquad x_i\in D_{p,q}^h.
\]
The selected coefficients may change independently at every node.
Scale row $i$ by $h^2/a_i>0$ and put
\[
r_i=\frac{h^2}{a_i},\qquad
\theta_i=\frac{hb_i}{2a_i},\qquad
g_j^o=\frac{h^2d_{2j+1}}{a_{2j+1}},\qquad
g_k^e=\frac{h^2d_{2k}}{a_{2k}}.
\]
With $v_k=U_{2k}$, $w_j=U_{2j+1}$, the odd equation contains only
its two adjacent even nodes, so it can be solved explicitly:
\begin{equation}
w_j=\frac{g_j^o+(1+\theta_{2j+1})v_j
+(1-\theta_{2j+1})v_{j+1}}{2+r_{2j+1}}.
\label{eq:q2-odd-recovery}
\end{equation}
The even equation before elimination is
\[
\tfrac14v_{k-1}-(2+\theta_{2k})w_{k-1}
+(\tfrac72+r_{2k})v_k
-(2-\theta_{2k})w_k+\tfrac14v_{k+1}=g_k^e.
\]
Substitute \eqref{eq:q2-odd-recovery}. Each formerly positive
$1/4$ coupling is offset by a path through the neighbouring odd node.
Collecting terms gives the tridiagonal Schur complement
\begin{equation}
-\eta_k^-v_{k-1}+(\eta_k^-+\eta_k^++\sigma_k)v_k
-\eta_k^+v_{k+1}=G_k,\qquad p<k<q,
\end{equation}
where
\[
\eta_k^-=
\frac{(2+\theta_{2k})(1+\theta_{2k-1})}{2+r_{2k-1}}-\frac14,
\qquad
\eta_k^+=
\frac{(2-\theta_{2k})(1-\theta_{2k+1})}{2+r_{2k+1}}-\frac14,
\]
\[
\sigma_k=r_{2k}
+\frac{(2+\theta_{2k})r_{2k-1}}{2+r_{2k-1}}
+\frac{(2-\theta_{2k})r_{2k+1}}{2+r_{2k+1}},
\]
\[
G_k=g_k^e
+\frac{(2+\theta_{2k})g_{k-1}^o}{2+r_{2k-1}}
+\frac{(2-\theta_{2k})g_k^o}{2+r_{2k+1}}.
\]
The diagonal identity follows by combining each eliminated-node
contribution with its adjacent off-diagonal coefficient; the remaining
positive reaction terms are exactly $\sigma_k$.

Put $\vartheta=hB/(2a_0)$. The mesh condition gives
$\vartheta<1/2$, $|\theta_i|\le\vartheta$ and $r_i\le h^2/a_0$.
In particular, all coefficients in the odd-node recovery and in
$G_k$ are nonnegative. More decisively,
\[
\begin{aligned}
4(2-\vartheta)(1-\vartheta)-(2+h^2/a_0)
&=6-12\vartheta+4\vartheta^2-h^2/a_0\\
&=\frac{6a_0-6hB-h^2}{a_0}+4\vartheta^2>0.
\end{aligned}
\]
Consequently
\[
\eta_k^\pm\ge
\frac{(2-\vartheta)(1-\vartheta)}{2+h^2/a_0}-\frac14>0,
\qquad \sigma_k>0.
\]
This calculation explains the sufficient mesh restriction. Only the
lower diffusion bound enters; arbitrarily large diffusion contrast
does not obstruct this one-dimensional elimination.

With $v_p,v_q$ prescribed, the reduced matrix has positive diagonal,
nonpositive off-diagonal entries and strict row diagonal dominance.
It is nonsingular: applying the equation at a component of maximum
absolute value rules out a nonzero kernel. It also has a nonnegative
inverse. Indeed, if its right-hand side is nonnegative and an interior
minimum $v_k<0$ existed with zero endpoint values, then
\[
\eta_k^-(v_k-v_{k-1})+\eta_k^+(v_k-v_{k+1})+\sigma_kv_k<0,
\]
a contradiction. Nonnegative endpoint values enter the reduced
right-hand side with coefficients $\eta_k^\pm>0$, so the even-node
solution is nonnegative whenever $d$ and the endpoint values are
nonnegative. Formula \eqref{eq:q2-odd-recovery} preserves this property
at odd nodes. If $q=p+1$, there are no interior even nodes, and the
same conclusion follows directly from the single odd equation.

Because $d$ was arbitrary and all row scalings are positive, this
proves $(A_D^\pi)^{-1}\ge0$ and
$-(A_D^\pi)^{-1}C_D^\pi\ge0$.
\end{proof}


{\small
\begin{thebibliography}{99}
	
	         \bibitem{AchdouFalcone2012}
		Y.~Achdou and M.~Falcone,
		\newblock A semi-Lagrangian scheme for mean curvature
		motion with nonlinear Neumann conditions,
		\newblock \emph{Interfaces and Free Boundaries}
		\textbf{14} (2012), 455--485.
		\newblock
		\href{https://doi.org/10.4171/IFB/288}
		{doi:10.4171/IFB/288}.
		
		\bibitem{BarlesDaLio2004}
		G.~Barles and F.~Da Lio,
		\newblock On the generalized Dirichlet problem for viscous
		Hamilton--Jacobi equations,
		\newblock \emph{Journal de Math\'ematiques Pures et Appliqu\'ees}
		\textbf{83} (2004), 53--75.
		
		\bibitem{BarlesJakobsen2002}
		G.~Barles and E.~R. Jakobsen,
		\newblock On the convergence rate of approximation schemes for
		Hamilton--Jacobi--Bellman equations,
		\newblock \emph{ESAIM: Mathematical Modelling and Numerical Analysis}
		\textbf{36} (2002), 33--54.
		
		\bibitem{BarlesSouganidis1991}
		G.~Barles and P.~E. Souganidis,
		\newblock Convergence of approximation schemes for fully nonlinear
		second order equations,
		\newblock \emph{Asymptotic Analysis}
		\textbf{4} (1991), 271--283.
		
		\bibitem{BokanowskiEtAl2015}
		O.~Bokanowski, M.~Falcone, R.~Ferretti, L.~Gr\"une,
		D.~Kalise, and H.~Zidani,
		\newblock Value iteration convergence of $\epsilon$-monotone schemes
		for stationary Hamilton--Jacobi equations,
		\newblock \emph{Discrete and Continuous Dynamical Systems}
		\textbf{35} (2015), 4041--4070.
		
		\bibitem{BokanowskiFalconeSahu2016}
		O.~Bokanowski, M.~Falcone, and S.~Sahu,
		\newblock An efficient filtered scheme for some first order
		time-dependent Hamilton--Jacobi equations,
		\newblock \emph{SIAM Journal on Scientific Computing}
		\textbf{38} (2016), A171--A195.
		
		\bibitem{BokanowskiPicarelliReisinger2018}
		O.~Bokanowski, A.~Picarelli, and C.~Reisinger,
		\newblock High-order filtered schemes for time-dependent second order
		HJB equations,
		\newblock \emph{ESAIM: Mathematical Modelling and Numerical Analysis}
		\textbf{52} (2018), 69--97.
		\newblock
		\href{https://doi.org/10.1051/m2an/2017039}
		{doi:10.1051/m2an/2017039}.
		
		\bibitem{CamilliFalcone1995}
		F.~Camilli and M.~Falcone,
		\newblock An approximation scheme for the optimal control
		of diffusion processes,
		\newblock \emph{RAIRO Mod\'el. Math. Anal. Num\'er.}
		\textbf{29} (1995), no.~1, 97--122.

		\bibitem{CarliniFalconeFerretti2010}
		E.~Carlini, M.~Falcone, and R.~Ferretti,
		\newblock Convergence of a large time-step scheme for mean curvature motion,
		\newblock \emph{Interfaces and Free Boundaries}
		\textbf{12} (2010), 409--441.
		\newblock
		\href{https://doi.org/10.4171/IFB/240}{doi:10.4171/IFB/240}.

		\bibitem{CarliniFerrettiPredaSemplice2025}
		E.~Carlini, R.~Ferretti, S.~Preda, and M.~Semplice,
		\newblock A CWENO large time-step scheme for
		Hamilton--Jacobi equations,
		\newblock \emph{Communications on Applied Mathematics and Computation}
		(2025).
		
		\bibitem{CrandallIshiiLions1992}
		M.~G. Crandall, H.~Ishii, and P.-L. Lions,
		\newblock User's guide to viscosity solutions of second order
		partial differential equations,
		\newblock \emph{Bulletin of the American Mathematical Society}
		\textbf{27} (1992), 1--67.
		
		\bibitem{FengLewis2021}
		X.~Feng and T.~Lewis,
		\newblock A narrow-stencil finite difference method for approximating
		viscosity solutions of Hamilton--Jacobi--Bellman equations,
		\newblock \emph{SIAM Journal on Numerical Analysis}
		\textbf{59} (2021), 886--924.
		
		\bibitem{FengLewisWard2022}
		X.~Feng, T.~Lewis, and K.~Ward,
		\newblock A narrow-stencil framework for convergent numerical
		approximations of fully nonlinear second order PDEs,
		\newblock \emph{Electronic Journal of Differential Equations,
			Conference}
		\textbf{26} (2022), 59--95.
		
		\bibitem{FroeseOberman2013}
		B.~D. Froese and A.~M. Oberman,
		\newblock Convergent filtered schemes for the Monge--Amp\`ere
		partial differential equation,
		\newblock \emph{SIAM Journal on Numerical Analysis}
		\textbf{51} (2013), 423--444.
		
		\bibitem{JensenSmears2013}
		M.~Jensen and I.~Smears,
		\newblock On the convergence of finite element methods for
		Hamilton--Jacobi--Bellman equations,
		\newblock \emph{SIAM Journal on Numerical Analysis}
		\textbf{51} (2013), 137--162.
		
		\bibitem{KuoTrudinger1992}
		H.-J.~Kuo and N.~S. Trudinger,
		\newblock Discrete methods for fully nonlinear elliptic equations,
		\newblock \emph{SIAM Journal on Numerical Analysis}
		\textbf{29} (1992), 123--135.
		\newblock
		\href{https://doi.org/10.1137/0729008}{doi:10.1137/0729008}.
		
		\bibitem{LiZhang2020}
		H.~Li and X.~Zhang,
		\newblock On the monotonicity and discrete maximum principle of the
		finite difference implementation of $C^0$-$Q^2$ finite element method,
		\newblock \emph{Numerische Mathematik}
		\textbf{145} (2020), 437--472.
		
		\bibitem{NochettoNtogkasZhang2019}
		R.~H. Nochetto, D.~Ntogkas, and W.~Zhang,
		\newblock Two-scale method for the Monge--Amp\`ere equation:
		convergence to the viscosity solution,
		\newblock \emph{Mathematics of Computation}
		\textbf{88} (2019), 637--664.
		
		\bibitem{Oberman2006}
		A.~M. Oberman,
		\newblock Convergent difference schemes for degenerate elliptic and
		parabolic equations: Hamilton--Jacobi equations and free boundary
		problems,
		\newblock \emph{SIAM Journal on Numerical Analysis}
		\textbf{44} (2006), 879--895.
		\newblock
		\href{https://doi.org/10.1137/S0036142903435235}
		{doi:10.1137/S0036142903435235}.
		
		\bibitem{PollockZhu2020Matrix}
		S.~N. Pollock and Y.~Zhu,
		\newblock A matrix analysis approach to discrete comparison principles
		for nonmonotone PDE,
		\newblock \emph{Numerical Algorithms}
		\textbf{83} (2020), 1007--1027.
		
		\bibitem{SalgadoZhang2019}
		A.~J. Salgado and W.~Zhang,
		\newblock Finite element approximation of the Isaacs equation,
		\newblock \emph{ESAIM: Mathematical Modelling and Numerical Analysis}
		\textbf{53} (2019), 351--374.
		\newblock
		\href{https://doi.org/10.1051/m2an/2018067}
		{doi:10.1051/m2an/2018067}.

		\bibitem{Warin2019}
		X.~Warin,
		\newblock Some non monotone schemes for
		Hamilton--Jacobi--Bellman equations,
		\newblock \emph{ESAIM: Proceedings and Surveys}
		\textbf{65} (2019), 476--497.
		
	\end{thebibliography}
}
\end{document}